\magnification\magstep1
\baselineskip=18pt
\let\a=\i
 
\def\li{\l_{\i}^n}
\def\i{\infty}
\def\l{\ell}
\def\t{\theta}
\def\la{\lambda}
\def\ra{\rightarrow}
\def\v{\vert}

\def\n{\noindent}
\def\Proof:{\n{\bf Proof:}}

\centerline{{\bf A simple proof of a
theorem of Jean Bourgain}\footnote*{Supported in part by
N.S.F. grant  DMS 9003550}} \vskip12pt \centerline {by  G.
Pisier} \vskip12pt

{\bf Abstract.}  
We give a simple proof of Bourgain's disc algebra version
of Grothendieck's theorem, i.e. that every operator on the
disc algebra with values in $L_1$
or $L_2$ is 2-absolutely summing and hence extends to an
operator defined on the whole of $C$. This implies
Bourgain's result that $L_1/H^1$ is of cotype 2. We
also prove more generally that   $L_r/H^r$  is of cotype
2 for $0<r< 1$.
\vfill\eject

In this note, we give a very simple proof (compared to
the preceding known proofs) of Bourgain's version of
Grothendieck's theorem for the disc algebra. As far as we
know, the currently known proofs are essentially the
original one in [B1], the simpler one in [BD], and several
new proofs given recently by Kisliakov, in [K1,K2]. 

We first recall the definition of a $q$-absolutely summing
(in short $q$-summing) operator for $1\leq q < \i$. Let
$u:X\rightarrow Y$ be an operator between two Banach
spaces. We say that $u$ is $q$-summing if there is a
constant $C$ such that for all finite sequences
$x_1,x_2,..,x_n$ in $X$, we have
$$(\sum\|u(x_i)\|^q)^{1/q} \leq C
\sup\{(\sum|x^*(x_i)|^q)^{1/q} |\  x^*\in X^* ,\  \|x^*\|
\leq 1 \}.  $$ We denote by $\pi_q(u)$ the smallest
possible constant $C$.  Let us denote by $A$ the disc
algebra. Then if $u:A\rightarrow Y$ is $q$-summing, by
Pietsch's factorisation theorem, there is a probability
measure $\lambda$ on the unit circle such that
$$\forall f \in A \qquad \|u(f)\| \leq
\pi_q(u) (\int{|f|^q d\lambda})^{1/q}.$$
 We refer  e.g. to
[P1] for more information on this notion.

We will prove

\noindent {\bf Bourgain's Theorem:}
There is a constant $K$ such that  any bounded operator
$u:A\rightarrow \ell_2$  is 2-summing and
satisfies:
$$\pi_2(u) \leq K \|u\|.$$
Also, $u$ extends to a bounded operator
$\hat{u}:C(\bf{T})\rightarrow \ell_2$ such that
$$\|\hat{u}\|\leq K\|u\|.$$
Moreover, the same result holds for all operators $u:A\rightarrow
Y$  if $Y=\ell_1$, or more generally, whenever $Y$ is a
Banach space of cotype 2.

Let us recall here the definitions of the $K_t$ and $J_t$
functionals which are fundamental in the real
interpolation method.
Let $A_0,A_1$ be a compatible couple of Banach (or
quasi-Banach) spaces. For all $x\in A_0+A_1$ and for all
$t>0$, we let $$K_t(x,A_0,A_1)=
\inf \big({\|x_0\|_{A_0}+t\|x_1\|_{A_1}\ | \
x=x_0+x_1,x_0\in A_0,x_1\in A_1}).$$
For all $x\in A_0\cap A_1$ and for all $t>0$, we let
$$J_t(x,A_0,A_1)=\max(\|x_0\|_{A_0},t\|x_1\|_{A_1}).$$
Recall that the (real interpolation) space
$(A_0,A_1)_{\theta,p}$ is defined as the space of all $x$
in $A_0+A_1$ such that $\|x\|_{\theta,p} <\infty$ where
 $$\|x\|_{\theta,p} =(\int{(t^{-\theta}K_t(x,A_0,A_1))^{p}
dt/t})^{1/p} .$$
We refer to [BeL] for more details.

Let {\bf T} be the circle group equipped with its
normalized Haar measure $m$. Let $1\le p\le \i$. When $B$
is a Banach space we denote by $L_p(B)$ the  usual space of
Bochner-$p$-integrable $B$-valued functions on $({\bf 
T},m)$, so that when $p<\i$, $L_p \otimes B$ is dense in
$L_p(B)$. We denote by ${H}^p(B)$ the subspace of
$L_p(B)$ formed by all the functions $f$ such that their
Fourier transform vanishes on the negative integers. When
$B$ is one dimensional, we write ${H}^p$ instead of
${H}^p(B)$. When $0<p<1$, we define ${H}^p$ as the closure
in $L_p$ of the linear span of the functions $\{e^{int}\ 
\v\  n\ge 0\}$. We refer to [G,GR] for basic information on
$H^p$-spaces.

 The next proposition although very
simple is the key new ingredient in our proof. We
refer to [P2] for more applications of the same idea to the
interpolation spaces between $H^p$ spaces.

\noindent {\bf Proposition 1:}\quad Let
$1\leq p\leq q\leq \infty$. Consider a compatible couple 
of Banach spaces $(A_0,A_1)$,  
the following are equivalent: 

$(i)$ There is a constant C such that
$$ \forall f \in H^p(A_0)+H^q(A_1), \quad\forall t>0,\quad
K_t(f,H^p(A_0),H^q(A_1)) \leq C
K_t(f,L^p(A_0),L^q(A_1)).$$

$(ii)$ There is a constant C such that
$$ \forall f \in [L^p(A_0)/H^p(A_0)]\cap
[L^q(A_1)/H^q(A_1)],\  \forall t>0,\quad  \exists \hat{f}
\in L^p(A_0)\cap L^q(A_1)$$ representing the equivalence
class of $f$ and satisfying
$$J_t(\hat{f},L^p(A_0),L^q(A_1))
\leq C J_t(f,L^p(A_0)/{H^p(A_0)},L^q(A_1)/{H^q(A_1)}).$$

$(iii)$ There is a constant C such that
$$ \forall f \in {[L^p(A_0)/ H^p(A_0)]}\cap
{[L^q(A_1)/H^q(A_1)]},  \quad  \exists
\hat{f} \in L^p(A_0)\cap L^q(A_1)$$ representing the equivalence
class of $f$ and satisfying
$$\|\hat{f}\|_{L^p(A_0)}\leq
C\|f\|_{L^p(A_0)/H^p(A_0)}\quad and \quad
\|\hat{f}\|_{L^q(A_1)}\leq C\|f\|_{L^q(A_1)/H^q(A_1)}.$$

In the above statement we regard the spaces 
$L^p (A_0)/H^p (A_0)$ and $L^q(A_1)/H^q(A_1)$ as included
via the Fourier transform $f\ra
(\hat{f}(-1),\hat{f}(-2),\hat{f}(-3),...)$ in the space of
all sequences in $A_0 + A_1$. In this way, we may view
these quotient spaces as forming a  compatible couple for
interpolation. (For the subspaces 
$H^p(A_0),H^q(A_1)$, there is no problem, we may clearly
consider them as a compatible couple in the obvious
way.)

 \Proof: For brievity, we will denote simply
$L^p/H^p(A_0)$ instead of $L^p(A_0)/H^p(A_0)$, we will
also write $L^p$,$H^p$,..instead of $L^p(A_0)$ ,
$H^p(A_0)$...no confusion should arise. The proof is
routine. We only indicate the argument for $(i)
\Rightarrow (ii) \Rightarrow (iii)$ which is the one we
use below.

Assume $(i)$.
Let $f$ be as above such that
 $J_t(f,L^p/H^p(A_0),L^q/H^q(A_1))<1$. Then let $g_p \in
L^p(A_0)$\quad  and $g_q \in L^q(A_1)$ be such that 
$$\|g_p\|_{L^p} <1, \quad \|g_q\|_{L^q} <t^{-1}, \quad 
f=g_p+H^p(A_0) ,\quad f=g_q+H^q(A_1).$$
Therefore, $g_p-g_q$ must be in $H^p+H^q$ and 
 $$K_t(g_p-g_q,L^p(A_0),L^q(A_1)) \leq
\|g_p\|_{L^p}+t\|g_q\|_{L^q}<2.$$
By $(i)$, we have $K_t(g_p-g_q,H^p,H^q)<2C'$, hence there
are $h_p\in H^p(A_0)$ and $h_q\in H^q(A_1)$ such that
$g_p-g_q=h_p-h_q$ and $\|h_p\|_{H^p}+t\|h_q\|_{H^q}<2C'$.
Now if we let $\hat f =g_p-h_p=g_q-h_q$,
then we find that $\hat{f} \in L^p(A_0)\cap
L^q(A_1)$, $f=\hat f +H^p(A_0)$ in the space
$L^p/H^p(A_0)$ and  $f=\hat f +H^q(A_1)$ in the space
$L^q/H^q(A_1)$ and moreover $$J_t(\hat f,L^p,L^q) \leq
\max(\|\hat f\|_{L^p},t\|\hat f\|_{L^q})\leq 1+2C'.$$ By
homogeneity this completes the proof of $(i)\Rightarrow
(ii)$ with $C\leq 1+2C'$. To check $(ii)\Rightarrow
(iii)$, simply write $(ii)$ with  $t=
(\|f\|_{L^p/H^p(A_0)}).(\|f\|_{L^q/H^q(A_1)})^{-1} .$

\noindent {\bf Remark :}
 It is well known that the Hilbert transform is a bounded
operator on all the (so-called mixed norm) spaces of the
form $L^p(\ell^q)$ for all $1<p,q<\infty$. (Apparently
this goes back to [BB]. We refer to [GR] for more
information and references). Therefore, the orthogonal
projection from $L^2(\ell^2)$ onto $H^2(\ell^2)$ is bounded
$\underline{simultaneously}$ on all the spaces
$L^p(\ell^q)$ for  $1<p,q<\infty$. It follows immediately
that if $1<p_0,p_1,q_0,q_1<\infty$ there is a constant $C$
such that  $ \forall f \in H^{p_0}(\ell
^{q_0})+H^{p_1}(\ell ^{q_1}), \quad\forall t>0,$ $$\quad
K_t(f,H^{p_0}(\ell ^{q_0}),H^{p_1}(\ell ^{q_1})) \leq C'
K_t(f,L^{p_0}(\ell ^{q_0}),L^{p_1}(\ell ^{q_1})).$$

\noindent {\bf Proposition 2:}
There is a constant $C'$ such that for all $t>0$ and all
$f\in H^1(\ell_1)+H^1(\ell_2)$ we have 
$$K_t(f,H^1(\ell_1),H^1(\ell_2)) \leq
CK_t(f,L^1(\ell_1),L^1(\ell_2)).$$

For the proof of Prop. 2, we will use the

\noindent {\bf Sublemma:}
$$H^1(\ell_{4/3}) \subset (H^1(\ell_1),H^1(\ell_2))_{1/2,\infty},$$
and the inclusion is bounded with norm less than a
constant $K$.

\Proof: Take a function $f=(f_n)$ in the unit ball of
$H^1(\ell_{4/3})$ and factor it as $f=(g_nh_n)$ with
$g=(g_n)$ in the unit ball of $H^2(\ell_2)$ and $h=(h_n)$ in
the unit ball of $H^2(\ell_4)$. This is easy to do by
factoring out the Blaschke product of each component $f_n$
and raising the factor without zero to the appropriate
power. (More precisely, write $f_n=B_nF_n$ where $B_n$ 
is a Blaschke product
and where $F_n$ does not have zeros in $D$, let $F$ be
an outer function such that we have
$\v F\v =(\sum \v F_n\v^{4/3})^{3/4}$
on the unit circle, then let
$$g_n=B_n(F_n/F)^{2/3} F^{1/2} \hbox{ and }
h_n=(F_n/F)^{1/3} F^{1/2}.$$
This factorisation has the   properties claimed for $g$
and $h$.) \vskip1pt
\n Recall the  inclusion (which obviously follows 
from the above remark) $$H^2(\ell_2)=
(H^2(\ell_{4/3}),H^2(\ell_4))_{1/2}\subset
(H^2(\ell_{4/3}),H^2(\ell_4))_{1/2,\infty} .$$ Then, by
interpolation, since the operator of coordinatewise
multiplication by $h=(h_n)$ maps $H^2(\ell_{4/3})$ into
$H^1(\ell_1)$ and $H^2(\ell_4)$ into $H^1(\ell_2)$, we
obtain the announced inclusion. q.e.d.

\noindent{\bf Proof of Prop. 2:}
Consider $f=(f_n) \in H^1(\ell_1)+H^1(\ell_2) $ such that
$$ K_t(f,L^1(\ell_1),L^1(\ell_2))<1.\leqno(1)$$
By classical factorisation theory, each $f_n$ can be
factored as $f_n=B_nF_n$ where $B_n$ is a Blaschke product
and where $F_n$ does not have zeros in $D$  so that
the analytic function $(F_n)^p$ makes sense
for any $p>0$.  (Alternatively, we could use the
inner-outer factorisation instead.) 
Let us simply denote by $F^{1/2}$ the sequence of
analytic functions $F^{1/2}=({F_n}^{1/2})_{n\geq1}.$
Note that any assumption of the form (1) depends only on
the values of each $|f_n|$ on the boundary. Now, on the
boundary we have  $|f_n|^{1/2}=|F_n|^{1/2}$,  so that (1)
implies obviously $$
K_{t^{1/2}}(F^{1/2},L^2(\ell_2),L^2(\ell_4))<2^{1/2}.\leqno(2)$$
Therefore, by the above Remark, 
$$K_{t^{1/2}}(F^{1/2},H^2(\ell_2),H^2(\ell_4))<2^{1/2}C,$$
where $C$ is a numerical (absolute) constant. Hence, there
is a decomposition $F^{1/2}=g_0+g_1$ with
$$\|g_0\|_{H2(\ell^2)}+{t^{1/2}}\|g_1\|_{H2(\ell^4)}
<2^{1/2}C.\leqno(3)$$ Let us now return to $f=(f_n)=(B_n
({g_0}_n+{g_1}_n)^2)$. Let us simply denote by $g_0g_1$ the
sequence $({g_0}_n{g_1}_n)_{n\geq1}$,similarly, we also
denote by ${g_0}^2$ and ${g_1}^2$ the sequences of
squares. Observe that by (3) and  by H\"older, we have
$$\|g_0g_1\|_{H^1(\ell_{4/3})} < 2C^2 t^{-1/2},$$ which
implies by the sublemma
$$t^{-1/2}K_t(g_0g_1,H^1(\ell_1),H^1(\ell_2)) < 2C^{2}K
t^{-1/2}.$$ After simplification
$$K_t(g_0g_1,H^1(\ell_1),H^1(\ell_2)) < 2C^{2}K.$$ On the
other hand we have clearly by (3)
 $$K_t({g_0}^2+{g_1}^2,H^1(\ell_1),H^1(\ell_2)) \leq
2C^2+2C^2=4C^2.$$

Therefore, we conclude by the triangle inequality (and
the fact that Blaschke products are of unit norm in
$H^\infty$)) $$K_t(f,H^1(\ell_1),H^1(\ell_2)) \leq
K_t({g_0}^2+{g_1}^2,H^1(\ell_1),H^1(\ell_2))
+K_t(2g_0g_1,H^1(\ell_1),H^1(\ell_2))$$
$$ \leq 4C^2+4C^2K.$$
By homogeneity, this completes the proof. q.e.d.

\def\l{\ell}
\def\t{\theta}
\def\la{\lambda}

\noindent {\bf Corollary:}
There is a constant $C$ such
that for all $1<p<2$ and all $f\in L^1/H^1(\l_p)$ we have

$$ \|f\|_{ L^1/H^1(\l_p)} \leq C\|f\|_{
L^1/H^1(\l_2)}^{\t}.\|f\|_{
L^1/H^1(\l_1)}^{1-\t}\leqno(4)$$ where $1/p=\t/2+(1-\t)/1$.

\Proof: By Prop. 2 and Prop. 1, there is a constant $C$
such that every $f\in L^1/H^1(\l_1)$ admits a
lifting  $\hat{f}\in L^1(\l_1)$ such that we have
$\underline {simultaneously}$
$$\|\hat{f}\|_{L^1(\l_1)} \leq C\|f\|_{L^1/H^1(\l_1)}$$
$$\|\hat{f}\|_{L^1(\l_2)} \leq C\|f\|_{L^1/H^1(\l_2)}.$$
Then (4) is an immediate consequence of H\"older's
inequality. q.e.d.

The preceding corollary implies immediately

\noindent {\bf Proposition 3:}
There is a constant $C$ such that, for all Banach spaces
$Y$,  for all $2<q<\infty$ and all 2-summing operator
$u:A\rightarrow Y$,  we have
$$\pi_q(u)\leq C\pi_2(u)^{\t} .\|u\|^{1-\t},\leqno(5)$$
where $1/q=\t/2+(1-\t)/\infty.$

\Proof: We first claim that for any $n>1$ and for any
$x_1,x_2,$...,$x_n$ in $A$, we have 
$$\sum_1^n \|u(x_i)\| \leq \la
\|(\sum|x_i|^q)^{1/q}\|_\infty,$$
where $\la\leq Cn^{1/q'}$.
Indeed, let us denote by $\la(q,n)$ the best constant in
this inequality. Assume w.l.o.g. that $u$ is the adjoint
of an operator $v: Y^*\rightarrow L^1/H^1 $. Let $p=q'$.
By duality, we find  $$\la(q,n)=\sup\{\|(v(y_i))\|_{
L^1/H^1(\l^n_p)}\}$$ where the sup runs over all n-tuples
$(y_i)$ in $Y$ such that $\sup\|y_i\| \leq 1$.
Therefore, (4) immediately yields
$\la(q,n)\leq C\la(2,n)^\t {\la(\infty,n)}^{1-\t}\leq
C(n^{1/2}\pi_2(u))^\t (n\|u\|)^{1-\t}$, hence $$ 
\la(q,n)\leq Cn^{1/q'}\pi_2(u)^\t \|u\|^{1-\t}.\leqno(6)$$
For simplicity, let $B=C\pi_2(u)^\t \|u\|^{1-\t}.$ By
(6),  we have for any $x_1,x_2,$...,$x_n$ in $A$,
$$ n^{-1/q'}\sum_1^n \|u(x_i)\| \leq
B\|(\sum|x_i|^q)^{1/q}\|_\infty.\leqno(7)$$
Now let us rewrite (7) for a sequence composed of
  $x_1/(k_1)^{1/q}$ repeated
$k_1$-times,$x_2/(k_2)^{1/q}$ repeated $k_2$-times,
etc.. We obtain $$(\sum k_i)^{-1/q'}\sum
{k_i}^{1/q'} \|u(x_i)\| \leq
B\|(\sum|x_i|^q)^{1/q}\|_\infty.$$ Clearly, since the
sequences of the form $((\sum k_i)^{-1} k_i)$ are
obviously dense in the set of all sequences $(\alpha_i)$
such that $\sum\alpha_i=1$, we obtain $$\sum
(\alpha_i)^{1/q'}\|u(x_i)\|\leq
B\|(\sum|x_i|^q)^{1/q}\|_\infty.$$ Taking the supremum over all
such $(\alpha_i)$, we finally obtain the announced result
(5). q.e.d.

We now recall a classical inequality due to Khintchine,
concerning the Rademacher functions $r_1,r_2,...,r_n,..$
defined on the Lebesgue interval. 
For every $q>2$, there is a constant $B_q$ such that
for all finite sequences of scalars $(\alpha_i)$, we have
$$(\int{|\sum\alpha_i r_i|^{q} dt})^{1/q} \leq B_q (\sum
|\alpha_{i}|^{2})^{1/2}.$$

The following is a  known result of Maurey [M].

\noindent {\bf Proposition 4:}
Let $X$ be any Banach space.
Let $Y$ be a Banach space of cotype 2, i.e. such that
there is a constant $C_2$ satisfying, for all $n$ and for
all $n$-tuples $y_1,y_2,..,y_n$ in $Y$,
$$(\sum \|y_i\|^2 )^{1/2}\leq {C_2}
 (\int{\|\sum r_i y_i \|^2 dt})^{1/2} .$$
Then, for every $q>2$, every
$q-$summing operator $u:X\rightarrow Y$ is actually
2-summing, and moreover
$$\pi_2(u) \leq B_q C_2 \pi_q(u).$$

\Proof: Let $x_1,x_2,..,x_n$ be a finite subset of $X$
such that $\sum|x^*(x_i)|^2 \leq 1 $ for all $x^*$ in the
unit ball of $X^*$. Then, by the above Khintchine
Inequality, we have for all $x^*$ in the
unit ball of $X^*$, $$(\int{|\sum r_{i}\  x^{*}(x_i)|^q dt}
)^{1/q} \leq B_q .$$
Hence, by the definition of $\pi_q(u)$ (note that the
integral below is actually an average over $2^n$ choices of
signs), $$ (\int{\|\sum r_i  u(x_i) \|^q dt})^{1/q} \leq 
\pi_q(u) B_q .$$
Hence, by the definition of the cotype 2,
$$(\sum \|u(x_i)\|^2)^{1/2} \leq C_2 B_q \pi_q(u) .$$
By homogeneity, this proves Proposition 4. q.e.d.

We can now complete the

\noindent {\bf Proof of Bourgain's Theorem:}
We use the same general line of attack as
Bourgain. This approach is based on an
extrapolation trick which originates in the work of Maurey
[M] and has been used several times before Bourgain's
work (especially by the author) to prove various extensions
of Grothendieck's Theorem. (The latter theorem corresponds
to the case $A=C$, $Y=\ell_1$ in the above statement, see
[P1] .) In this approach, the crucial point
reduces to showing (5). 
Indeed, assuming (5), it is easy to conclude:
By Prop. 4, we have $\pi_2(u)
\leq C_2 B_q \pi_q(u)$, hence by
(4),  $\pi_2(u)\leq C C_2 B_q\pi_2(u)^\t .\|u\|^{1-\t}$
hence if we assume a priori that $\pi_2(u)$
is finite, we obtain $$\pi_2(u)\leq (CC_2 B_q )^{1/{1-\t}}
\|u\| ,\leqno(8)$$ which establishes the announced result
in the case of a 2-summing  operator. Hence in particular
(8) holds if $u$ is of finite rank, so that we can easily
conclude, since $A$ has the Metric Approximation Property,
that it actually holds for arbitrary operators.
Finally, the last assertion follows from a well
known factorisation  property of 2-summing operators, due
to Pietsch, cf. e.g. [P1], chapter 1. q.e.d.

\noindent {\bf Remark 5:} There is also a
slightly different way to prove (8). One can use a simple
interpolation argument to prove  that for any $n>1$ and for
any $x_1,x_2,$...,$x_n$in $A$, we have
$$\|(u(x_i))\|_{\l_{q,\i}(Y)} \leq B \|(x_i)\|_{L_\i
(\l_{q,\i})}, $$
where $B$ is as above. Then, we may apply this replacing
 $x_1,x_2,$...,$x_n$ by the $2^n$-tuple 
formed by the $2^n$ choices of signs $\sum r_i(t)  x_i $.
After normalisation by a factor $2^{-n/q}$, we obtain
$$\|(\sum r_{i} u(x_i))\|_{L_{q,\i}(dt;Y)} \leq B
\|\sum r_i x_i\|_{L_\i (L_{q,\i}(dt))}. $$
But then, we observe that Khintchine's inequality implies
a fortiori the equivalence of $\|\sum r_i x_i\|_{L_\i
(L_{q,\i}(dt))} $ with $\|(x_i)\|_{L_\i (\l _2)}$.
This immediately leads to (8) by the same argument as
above.

\noindent {\bf Remark 6:} As is well known, it follows
by standard arguments from Bourgain's theorem as stated
above that $L_1/H^1$ is a cotype 2 space. This can be
derived as in [B1] from a result of Wojtaszczyk which
ensures that $L_1/H^1$ is isomorphic to $L_1/H^1 (\l_1)$.
Alternately, if one wishes to avoid the use of the latter
result, one can observe that our proof of Bourgain's
Theorem is valid with essentially the same proof with
$L_1/H^1 (\l_1)$ instead of $L_1/H^1$.

Actually, we can generalize Bourgain's theorem as follows:

\noindent {\bf Theorem 7:} Let $0<r<1$. Then  any operator 
 $u:c_0 \ra L_r/H^r $
is 2-summing. Moreover, there is a constant $C_r$ such that
every operator $u:c_0 \ra L_r/H^r $ is 2-summing and
satisfies $\pi_2(u) \leq C_r \|u\|$. Finally, $L_r/H^r$ is
of cotype 2 .

\Proof: We only sketch the argument. (It might very well be
that this result follows from the other proofs, however
it seems to have passed unnoticed so far.)
Consider an operator $u:\l_{\i}^n \ra L_r/H^r $. We will
show that there is a constant $C_r$ independent of $n$
such that, $$\forall m\quad \forall x_1,x_2,...,x_m \in \li
,\quad \|(u(x_i))\|_{L_r/H^r(\l_2^m)} \leq C_r \|u\|
\|(x_i)\|_{\li (\l_2^m)}.\leqno (9) $$ 
We argue similarly as above, but in a dual setting. Let
$r \leq p\leq \i$. We denote by $C_p(u)$ the smallest
constant $C$ such that  $$\forall m \quad\forall x_1,x_2,...,x_m
\in \li ,\quad \|(u(x_i))\|_{L_r/H^r(\l_p^m)} \leq C 
\|(x_i)\|_{ \l_p^m(\li)}. $$ 
Obviously, we have $C_r(u)= \|u\|$. Choose $p$ such that
$r<1<p<2.$ Let $\t$ be such that ${1/p}={(1-\t)}/r
+\t/2.$  A simple adaptation of Proposition 1 and 2
yields a simultaneous "good" lifting for the couple
$L_r/H^r(\l_r),L_r/H^r(\l_2)$, and the corresponding
extension of (4). It follows that we have for some constant
$C'$ (independent of $m$) $$\|(u(x_i))\|_{L_r/H^r(\l_p^m)}
\leq C' {C_2(u)}^{\t}\|u\|^{1-\t} m^{1/p} \sup
\|x_i\|_{\li}. $$ As a consequence, if $B'=  C'
{C_2(u)}^{\t}\|u\|^{1-\t}$, we have 

$$\|(u(x_i))\|_{L_r/H^r(\l_p^m)} \leq B'
\|(x_i)\|_{\l_{p,r}^{m} (\li)}. \leqno (10)$$ 
It is easy to check that for some constant $C''$
(independent of $m$ or $n$) we have,
$$\|(x_i)\|_{\l_{p,r}^{m} (\li)} \leq C''m^{1/p-1/2}
\|(x_i)\|_{\l_2^{m} (\li)} $$ so that (10) gives after
normalisation (here $L_p^{m}$ denotes the $L_p$-space
relative to $\{1,2...,m\}$ 
equipped with the uniform probability measure)
$$\|(u(x_i))\|_{L_r/H^r(L_p^m)} \leq B'C''
\|(x_i)\|_{L_2^{m} (\li)}. \leqno (11)$$ 
Let $K=B'C''$.
 Now we take $m=2^k$, we replace $(x_i)$ by the $2^k$
"choices of signs" $x_t = \sum_1^k r_i(t) x_i $ and
we use the dualisation of Khintchine's inequality in $L_p$
 which says that, if $p>1$,  the quotient of $L_p$ by the
orthogonal of the Rademacher functions can be identified
with $\l_2$. If we simply denote by $Q(n)$ the quotient
space of  $L_2(\li)$ by the subspace of all functions
"orthogonal" to the Rademacher functions ( i.e. which have
a zero integral against any Rademacher function), we
can deduce from (11)
$$\|(u(x_i))\|_{L_r/H^r(\l_2^k)} \leq K
\|(x_i)\|_{Q(n)}. \leqno (12)$$ 
But on the other hand by a known reformulation of
Grothendieck's theorem (see [P1], corollary 6.7, p. 77), we
have  $$\|(x_i)\|_{Q(n)} \leq K' \|(x_i)\|_{\li(\l_2^k)}
\leqno(13)$$ where $K'$ is a numerical constant.
Therefore, (12) implies
$$C_2(u) \leq K K'  .$$
Recalling the value of $K$ and $B'$, we conclude that
$$C_2(u) \leq K'C'' C'
{C_2(u)}^{\t} \|u\|^{1-\t} ,$$
so that we again conclude by "extrapolation"
that $C_2(u) \leq K'' \|u\| $ for some constant $K''$
depending only on $r$. Combining (12) and (13) with this
last estimate, we obtain the announced result (9) with
$C_r=K C' C''$. Since there is obviously a norm one
inclusion of ${L_r/H^r(\l_2^m)}$ into $\l_2^m(L_r/H^r)$,
we have $\pi_2(u) \leq C_2(u) \leq C_r \|u\|$, and this
completes the proof for $X=\li$, (with a constant $C_r$
bounded independently of $n$). By density, this is enough
to prove the case of an operator defined on $c_0$. Finally,
the cotype 2 property can be proved as indicated in Remark
6, by observing that the first part of Theorem 7 remains
valid with ${L_r/H^r(l_r)}$ (or equivalently
${l_r(L_r/H^r)}$ )  in the place of $L_r/H^r$. We then
follow a standard argument, given elements $x_1,x_2,..,x_n
\in L_r/H^r$, we consider the operator $u:\li \ra
L_r(L_r/H^r) $ defined by $u(a_1,a_2,...,a_n) =\sum a_i
r_i x_i $, where 
 $r_1,r_2,...,r_n$ are the Rademacher functions as
before. We have
$$(\sum\|x_i\|^2)^{1/2} \leq \pi_2(u)\leq C_r \|u\|
,\leqno(14)$$ but it is well known that there is a constant
$B_r$, depending only on $r$ such that 
$$ \|u\| \leq B_r\|\sum r_i(t) x_i\|_{L_r(dt;L_r/H^r)}    ,
$$
therefore, (14) implies that $L_r/H^r$ is a cotype 2
space. q.e.d.

\vfill\eject
 
\noindent{\bf Final Remarks:} 

1) As a
corollary, one obtains that every rank $n$ operator on $A$
extends to the whole of $C({\bf T})$ with norm at most $C
Log n$ for some constant $C$. This follows from Bourgain's
theorem and a previous result of Mityagin and
Pelczy\'nski, see [B1] for the deduction.

2) The preceding argument shows that $$H^\infty
(\l_{p,\infty}) = (H^\infty
(\l_{1}),H^\infty
(\l_{\infty})) _{\t,\infty} \leqno (15)$$
where $1/p=1-\t ,0<\t<1.$
But this kind of result is not really new. It can be
derived from the remarks on interpolation spaces included
in [B1] using a  rather simple factorisation argument,
such as for instance the one used for Theorem 2.7 in [HP].
More results along this line have been obtained by Xu
[X]. In [P2], we will give a more systematic treatment of
 results such as (15), in more general cases for the real
interpolation method with arbitrary parameters.

3) We should mention that while Kisliakov's recent proof
of Bourgain's theorem seems more complicated than the
above, it also yields more information (on the so-called
$(p,q)$-summing operators) which do not follow from our
approach, cf.[K2,K3]. Moreover, although the above argument
applies also for an operator defined on $H^\infty$ and with
values in a cotype 2 space $Y$ with the Bounded
Approximation Property, it is a well known drawback of the
"extrapolation method" that it does not apply to the case
of a linear operator from $H^\infty$ into its dual,
although that case was settled in [B2].

\centerline {\bf References}\vskip6pt

\item {[BeL]} J. Bergh and J. L\"ofstr\"om, Interpolation
spaces, An introduction, Springer Verlag 1976.

\item {[BB]} R. P. Boas and S. Bochner, On a theorem of
Marcel Riesz for Fourier series.J.London Math. Soc.
14(1939) 62-73.
 
\item {[B1]} J. Bourgain, New Banach space properties of
the disc algebra and $H^\infty$, Acta Math. 152 (1984)
1-48.

\item {[B2]} J. Bourgain, Bilinear forms on $H^\infty$ and
bounded bianalytic functions. Trans. Amer. Math. Soc. 286
(1984) 313-337.

\item {[BD]} J. Bourgain and W. J. Davis, Martingale
transforms and complex uniform convexity. Trans. Amer.
Math. Soc. 294 (1986) 501-515.

\item {[G]} J.Garnett, Bounded Analytic Functions.
Academic Press 1981.

\item {[GR]} J. Garcia-Cuerva and J. L. Rubio de Francia.
Weighted norm inequalities and related topics. North
Holland, 1985.

\item {[HP]} U. Haagerup and G. Pisier, Factorization of
analytic functions with values in non-commutative
$L_1$-spaces and applications. Canadian J. Math. 41
(1989) 882-906.

\item {[K1]} S.V.Kisliakov, Truncating functions in
weighted $H^p$ and two theorems of J.Bourgain. Preprint, 
Uppsala University, May 1989.

\item {[K2]} S.V.Kisliakov, $(q,p)$-summing operators on
the disc algebra and a weighted estimate for certain
outer functions. Lomi preprint. Leningrad,1989.

\item {[K3]} S.V.Kisliakov, Extension of $(q,p)$-summing 
operators defined on
the disc algebra, with an appendix on Bourgain's analytic
projection. Preprint.1990.

\item {[M]} B. Maurey, Th\'eor\`emes de factorisation
pour les op\'erateurs \`a valeurs $L^p$. Ast\'erisque, Soc.
Math. France. 11 (1974)

\item {[P1]} G. Pisier,  Factorization of linear
operators and geometry of Banach spaces. Amer. Math. Soc.
CBMS 60. 1986.

\item {[P2]} G. Pisier, Interpolation between $H^p$ spaces
and non-commutative generalizations. To appear.

\item {[X]} Quanhua Xu, Real interpolation of some Banach
lattices valued Hardy spaces.Pre-print. Pub. IRMA, Lille.
1990.vol. 20 n$^\circ$ 8.
\vskip12pt

Texas A. and M. University

College Station, TX 77843, U.S.A.

and

Universit\'e Paris 6
 
Equipe d'Analyse, Bo\^{\a}te 186,

75230 Paris Cedex 05, France

 \bye